\documentclass[3p,10pt,a4paper,twoside,fleqn,sort&compress]{filomat-arXiv}

\usepackage{amssymb,amsmath,latexsym}
\usepackage[varg]{pxfonts}

\newtheorem{theorem}{Theorem}[section]
\newtheorem{statement}{Statement}[section]

\newtheorem{remark}[theorem]{Remark}

\newcommand{\RR}{{\mathbb R}}

\newcommand{\NN}{{\mathbb N}}

\begin{document}

\title{One method for proving  some  classes of exponential analytical inequalities}

\author[ETF]{Branko Male\v sevi\' c${}^{\ast}$}
\ead{branko.malesevic@etf.bg.ac.rs}
\author[ETF]{Tatjana Lutovac}
\ead{tatjana.lutovac@etf.bg.ac.rs}
\author[ETFFTN]{Bojan Banjac}
\ead{bojan.banjac@uns.ac.rs }

\address[ETF]{School of Electrical Engineering, University of Belgrade, Serbia}
\address[ETFFTN]{Faculty of Technical Sciences, University of Novi Sad, Serbia\\
                 PhD student at School of Electrical Engineering, University of Belgrade, Serbia}
\newcommand{\AuthorNames}{B. Male\v sevi\' c, T. Lutovac, B. Banjac}

\newcommand{\FilMSC}{Primary 26D07; Secondary 41A58, 26E05}
\newcommand{\FilKeywords}{Analytic inequalities, Series expansions, Real-analytic functions}
\newcommand{\FilCommunicated}{Miodrag Spalevi\' c}
\newcommand{\FilSupport}{The first author was supported in part by the  Serbian Ministry of Education,
Science and Technological Development, under projects ON 174032 and III 44006.
The second author was supported in part by the Serbian Ministry of Education, Science and Technological Development, under project  TR 32023.}
\newcommand{\FilCorespond}{${}^{*}$ }

\begin{abstract}
In this paper we propose a method for proving  some  exponential inequalities based on
power series expansion and  analysis of derivations of the corresponding functions.
Our approach provides a simple proof and  generates  a new class of appropriate inequalities,
as well as allows  direct establishment of the dependence between (the exponent of)  some functions that occur as  bounds of the approximation
and the interval in which the corresponding inequality holds true.
\end{abstract}

\maketitle
\makeatletter
\renewcommand\@makefnmark%
{\mbox{\textsuperscript{\normalfont\@thefnmark)}}}
\makeatother

\section{Introduction and preliminaries} 

Papers \cite{Telfor2015}-\cite{RLM_2018} have recently considered various methods for proving
mixed trigonometric polynomial inequalities of the form:
$$
f(x)=\sum_{i=1}^{n}{\alpha_i x^{p_i} \sin ^{q_i} x \cos ^{r_i} x}>0,
$$
where $p_i,q_i,r_i \in \NN_0$, $\alpha_i \in \RR \;\backslash\!\left\{0\right\}$ and $x\in(0,\frac{\pi}{2})$.
In monographs \cite{Mitrinovic1970} and \cite{Milovanovic_Rassias_2014} were stated various analytical inequalities
that can be reduced to mixed trigonometric polynomial inequalities. An algorithm that reduces proving of such inequalities
to proving of the corresponding  polynomial inequalities is developed in  \cite{Lutovac_Malesevic_Mortici_2017}. It is shown
that many open problems and various inequalities recently published in renowned journals can be proved using the proposed algorithm.

\smallskip
In this paper we consider some  exponential analytical inequalities whose proving can be reduced to analysis
of power series expansion of the corresponding functions, and  analysis of  the derivatives of the function.
We propose a method  for generating and proving  an appropriate class of exponential inequalities.

\smallskip
The starting point in our analysis  are the following two (recently published) results:

\break

\begin{statement}
\label{Statement-1}
$($ {\rm \cite{Bagul2017}}$\,,$ Theorem $1$ $)$
If $x \in \left(0,1\right)$ then
\begin{equation}
\label{inequ_1}
e^{-ax^2} < \cos x < e^{-x^2/2}
\end{equation}
where $a \approx 0.615626$.
\end{statement}
\begin{statement}
$($ {\rm \cite{Bagul2018}}$\,,$ Theorem $3$  $)\,$
If $x \in \left(0,1\right)$ then
\begin{equation}
\label{inequ_2}
e^{-bx^2} < \frac{x}{\tan x} < e^{-x^2/3}
\end{equation}
with the best possible constants $b \approx 0.443023$ and $1/3$.
\end{statement}

In this paper we present and provide proof for an  improvement  of the previous statements.

\smallskip
 Let us recall some of the well-known power series expansions that will be used in our proofs.

\smallskip
Based on  \cite{Gradshteyn-Ryzhik},  the following power series expansions hold:
\begin{equation}
\label{Series_ln_cos_x}
\log \cos x
=
-\sum\limits_{k=1}^{\infty}{\frac{2^{2k-1}(2^{2k}-1)|\mbox{\bf B}_{2k}|}{k(2k)!}x^{2k}}, \qquad (-\pi/2 < x < \pi/2),
\end{equation}
and
\begin{equation}
\label{Series_ln_tan_x_over_x}
\log \frac{x}{\tan x}
=
-\sum\limits_{k=1}^{\infty}{\frac{(2^{2k-1}-1)|\mbox{\bf B}_{2k}|}{k(2k)!}x^{2k}}, \qquad (0 < x < \pi/2),
\end{equation}
where $\mbox{\bf B}_{i}$ ($i \!\in\! \NN$) are {\sc Bernoulli}'s numbers.

\section{Main results}

\smallskip
First, we consider the relationship between the number of zeros and the number of (local) minimums of a real function,
as well as some properties of its derivatives.

\smallskip
We prove the following assertion:

\begin{theorem}
\label{nike_theorem}
Let $ f: (0, c) \longrightarrow \RR $ be  $m$ times differentiable function $(\;$for some $m \geq 2, \, m \in  \NN{\,)}$
which satisfies the following conditions:
\begin{itemize}
\item[{\rm (a)}]  $f^{(m)}(x) > 0\,$ for $x \in (0,c)$;
\item[{\rm (b)}] there is a right neighbourhood  of zero in which  the following inequalities hold true:
$$\, f(x)<0, \, f'(x)<0, \, \ldots, f^{(m-1)}(x)<0;$$
\item[{\rm (c)}] there is a left neighbourhood of $c$ in which the following inequalities hold true:
$$
\, f(x)>0, \, f'(x)>0, \, \ldots, f^{(m-1)}(x)>0.
$$
\end{itemize}
Then the function $f$ has exactly one zero $x_0 \in (0, c)$,
and  $f(x)<0$ for $x \in (0,x_0)$ and $f(x)>0$ for $x \in (x_0, c)$. Also, the function $f$ has
exactly one local minimum   in the   interval  (0, c) i.e. there is  exactly one point
$ t_0 \in \left(0, x_ {0}\right)\subset \left(0, c\right)$ such that  $f(t_0)<0$ is the minimal value of function $f(x)$ over the interval
$\left(0, x_0 \right)$ i.e. $  \left(0, c\right)$.
\end{theorem}
{\bf Proof.}
On the basis of condition (a), as $f^{(m)}(x) > 0$ for $x \in (0,c)$, it follows that $f^{(m-1)}(x)$ is a monotonically increasing function
for $x\in (0,c)$. Based on conditions (b) and (c), we conclude that there exists exactly one zero $x_{m-1} \in (0, c)$ of the function $f^{(m-1)}(x)$.
Next, we can conclude that function $f^{(m-2)}(x)$ is  monotonically decreasing  for  $x \in (0,x_{m-1})$ and monotonically increasing for
$x \in (x_{m-1},c)$. It is clear that the function $f^{(m-2)}(x)$ has exactly one minimum in the interval $(0, c)$ at point $x_{m-1}$,
and so $f^{(m-2)}(x_{m-1})<0$  holds. On the basis of condition (c), it follows that function $f^{(m-2)}(x)$ has exactly one root
$x_{m-2}$ in the interval $(0, c)$  ($x_{m-2}\in (x_{m-1}, c)$) and  $f^{(m-2)}(x)<0$  for $x \in (0, x_{m-2})$ as well as
$f^{(m-2)}(x)>0$  for $x \in (x_{m-2},c)$.

\noindent
By repeating the above procedure, we get the assertion given in the theorem.\hfill $\Box$

\medskip
\begin{remark}
The  Theorem~\ref{nike_theorem} is a natural extension of  Theorem~$3$, from {\rm \cite{RLM_2018}}, which was
 applied in {\rm \cite{RLM_2018}} to inequalities involving the sinc function.
\end{remark}

\subsection{Generalization of  Statement \mbox{1.1}}

\medskip
Consider inequalities (\ref{inequ_1}) for $x \in  (0, \pi/2)$ and  parameter $a \in \RR^+$.

\smallskip
It is enough to analyse the following equivalent inequalities:
\begin{equation}
\label{logaritamTheorem1}
-ax^2<\log \cos x< -x^2/2,
\end{equation}
for  $x \in (0,\pi/2)$  and parameter $a \in \RR^+$.

\smallskip
Firstly, we prove that the right-hand side of inequality  (\ref{logaritamTheorem1}) holds true for every $x \in (0, \pi/2)$.

\smallskip
The real analytical function
$$
f_1(x)
=
-\mbox{\scriptsize $\footnotesize \displaystyle\frac{1}{2}$}x^2 - \log \cos x
:
\left(0,\mbox{\scriptsize $\footnotesize \displaystyle\frac{\pi}{2}$}\right) \longrightarrow \RR,
$$
has the following power series expansion
$$
f_1(x)
=
\sum\limits_{k=2}^{\infty}{\frac{2^{2k-1}(2^{2k}-1)|\mbox{\bf B}_{2k}|}{k(2k)!}x^{2k}},
$$
for  $x \in (0,\pi/2)$, from which it follows that $f_1(x) > 0$ for $x \in (0,\pi/2)$.

\medskip
Consider now the left-hand side of inequality (\ref{logaritamTheorem1}). The corresponding real analytical function
$$
f_2(x)
=
-ax^2 - \log \cos x
:
\left(0,\mbox{\scriptsize $\footnotesize \displaystyle\frac{\pi}{2}$}\right) \longrightarrow \RR
$$
has the following power series expansion:
\begin{equation}
\label{power-series-f_2}
f_2(x)
=
\left(\frac{1}{2}-a\right)x^2
+
\sum\limits_{k=2}^{\infty}{\frac{2^{2k-1}(2^{2k}-1)|\mbox{\bf B}_{2k}|}{k(2k)!}x^{2k}}.
\end{equation}
Based on the above power series expansion, it is not hard to show that for $a>1/2$  and $m=3$  the conditions of Theorem~\ref{nike_theorem}
are fulfilled. Therefore, in the interval $\left(0, \mbox{\scriptsize $ \footnotesize \displaystyle \frac{\pi}{2} $} \right) $
the function $f_2(x)$ has exactly one zero $x_0$, and
there is exactly one point $ t_0 \in \left(0, x_ {0} \right) $ such that $f_2(t_0)<0$  is minimal value of the function $f_2(x)$
over the  interval $\left(0,x_ {0} \right) \subset \left(0, \pi/2 \right)$.

\smallskip
For every fixed  $x_0 \in \left(0,\mbox{\scriptsize $\footnotesize \displaystyle\frac{\pi}{2}$}\right)$ and
$$
a=-\frac{\log \cos x_0}{x_0^2},
$$
we have $f_2(x_0)=0$, i.e. $x_0$ is  the (unique) zero of the function $f_2(x)$, and for all $x\in (0, x_0)$ $\,f_2(x)<0$.
Further, based on (\ref{power-series-f_2}) and $f_2(x_0)=0$ we have:
$$
a
=
\mbox{\small $\displaystyle\frac{1}{2}$}
+
\mbox{$\displaystyle\sum$}_{k=2}^{\infty}{\mbox{\scriptsize $\footnotesize \displaystyle\frac{2^{2k-1}(2^{2k}-1)|\mbox{\bf B}_{2k}|}{k(2k)!}$}
x_0^{2k-2}}
>
\mbox{\small $\displaystyle\frac{1}{2}$}.
$$

\medskip
Overall, on the basis of the previous consideration, the following improvement of Statement 1.1  has been
proved:

\begin{theorem}
\label{Theorem-1}
Let  $x_0 \in \left(0, \pi/2\right)$. Then for every $x \in \left(0, x_0 \right) $
the following inequalities hold true
$$ 
e^{-ax^2} < \cos x < e^{-x^2/2}
$$ 
with the best possible constants $a=-\left(\log \cos x_0\right)\,\mbox{\large /}\,x_0^2\,$ and $\,1/2$.
\end{theorem}

\begin{remark}
For $x_0=1$ we get  Statement {\rm 1.1.} and constant  $a=-\log \cos 1 \approx 0.615626$. In that case, $t_0 \approx 0.736713$
and the minimum value $f_2 (t_0) \approx -0.0339435$ determines the absolute value of the largest error $\delta \approx 0.0339435$
in the approximation of the function $ \log\cos x $ with the square polynomial $-ax^2$ over interval $(0,1)$.
\end{remark}


\subsection{Generalization of Statement \mbox{1.2}}

\medskip
Consider  inequalities (\ref{inequ_2}) for $x \in  (0, \pi/2)$ and  parameter $b \in \RR^{+}$.

\smallskip
It is enough to analyse the following equivalent inequalities:
\begin{equation}
\label{logaritamTheorem2}
-bx^2<\log \frac{x}{\tan x} <-\frac{x^2}{3}
\end{equation}
for $x \in  (0, \pi/2)$ and  parameter $b \in \RR^{+}$.

\smallskip
For the right-hand side of inequality (\ref{logaritamTheorem2}) it is enough to observe the real analytical  function
$$
g_1(x) = -\mbox{\scriptsize $\footnotesize \displaystyle\frac{1}{3}$}x^2 - \log \frac{x}{\tan x}
:
\left(0,\mbox{\scriptsize $\footnotesize \displaystyle\frac{\pi}{2}$}\right) \longrightarrow \RR,
$$
and the corresponding  power series expansion:
$$
g_1(x)
=
\sum\limits_{k=2}^{\infty}{\frac{2^{2k}(2^{2k-1}-1)|\mbox{\bf B}_{2k}|}{k(2k)!}x^{2k}},
$$
for $x \in (0,\pi/2)$,  from which it follows that $g_1(x) > 0$ for $x \in (0,\pi/2)$.

\medskip
Consider now the left-hand side of inequality  (\ref{logaritamTheorem2}).
The corresponding real\\ analytical function
$$
g_2(x)
=
-bx^2 - \log \frac{x}{\tan x}
:
\left(0,\mbox{\scriptsize $\footnotesize \displaystyle\frac{\pi}{2}$}\right) \longrightarrow \RR
$$
has the following power series expansion:
\begin{equation}
\label{power-series-g_2}
g_2(x)
=
\left(\frac{1}{3}-b\right)x^2
+
\sum\limits_{k=2}^{\infty}{\frac{2^{2k}(2^{2k-1}-1)|\mbox{\bf B}_{2k}|}{k(2k)!}x^{2k}}.
\end{equation}
It is not hard to show that  for $b>1/3$  and $m=3$  the conditions of Theorem~\ref{nike_theorem} are fulfilled.
Therefore,
the function $g_2(x)$ has exactly one zero $x_0$   in the interval $\left(0, \mbox{\scriptsize $ \footnotesize \displaystyle \frac{\pi}{2} $} \right) $,
and there is exactly one point $t_0 \in \left(0, x_ {0} \right)$
such that $g_2(t_0)<0$  is the minimal value of the function $g_2(x)$ over the interval
$\left(0, x_ {0} \right)\subset \left(0, \pi/2 \right)$.

\smallskip
For every fixed  $x_0 \in \left(0,\mbox{\scriptsize $\footnotesize \displaystyle\frac{\pi}{2}$}\right)$ and
$$\displaystyle
b=-\frac{\log \frac{x_0}{\tan x_0}}{x_0^2},
$$
we have $g_2(x_0)=0$, i.e. $x_0$ is the (unique) zero of the function $g_2(x)$, and for all $x\in (0, x_0)$ $\,g_2(x)<0$.
Further, based on (\ref{power-series-g_2}) and $g_2(x_0)=0$ we have:
$$
b
=
\mbox{\small $\displaystyle\frac{1}{3}$}
+
\mbox{$\displaystyle\sum$}_{k=2}^{\infty}{\mbox{\scriptsize $\footnotesize \displaystyle\frac{2^{2k}(2^{2k-1}-1)|\mbox{\bf B}_{2k}|}{k(2k)!}$}
x_0^{2k-2}}
>
\mbox{\small $\displaystyle\frac{1}{3}$}.
$$

\medskip
Overall, on the basis of the previous consideration, the following improvement of Statement 1.2  has been proved:

\begin{theorem}
\label{Theorem-1}
Let  $x_0 \in \left(0, \pi/2\right)$. Then for every $x \in \left(0, x_0 \right) $
the following inequalities hold true
$$ 
e^{-bx^2} < \log \frac{x}{\tan x} < e^{-x^2/3}
$$ 
with the best possible constants $b=-\left(\log \mbox{\footnotesize $\displaystyle\frac{x_0}{\tan x_0}$}\right)\,\mbox{\large /}\,x_0^2\,$
and $\,1/3$.
\end{theorem}

\begin{remark}
For $x_0=1$ we get  Statement 1.2 and constant $b=\log \tan 1 \approx 0.443023$. In that case, $t_0 \approx 0.737815$ and the minimum value
$g_2 (t_0) \approx -0.0324168$ determines the absolute value of the largest error $\delta \approx 0.0324168$ in the approximation of the
function $\log \mbox{\footnotesize $\displaystyle\frac{x}{\tan x}$}$ with the square polynomial $-bx^2$ over interval $(0,1)$.
\end{remark}

\section{Conclusion}
In this paper we  established  some exponential analytic inequalities, using the results of Theorem~\ref{nike_theorem},
as well as the power series expansion of the corresponding functions. In particular,
we proved that the following inequalities hold true:
$$
\begin{array}{l}\displaystyle
e^{-ax^2} < \cos x < e^{-x^2/2}, \\ [1em]
e^{-bx^2} < \log \frac{x}{\tan x} < e^{-x^2/3}
\end{array}
$$
for every $x \in \left(0, \alpha \right),\, $ where $0< \alpha < \pi/2$,
with the best possible constants $a=-\left(\log \cos \alpha \right)\,\mbox{\large /}\,\alpha^2\,$ and $\,1/2$,
and  $b=-\left(\log \mbox{\footnotesize $\displaystyle\frac{\alpha}{\tan \alpha}$}\right)\,\mbox{\large /}\,\alpha^2\,$  and $\,1/3$, respectively.

Our approach provides a simple proof and  allows  direct establishment of the dependence between
the exponents  of  the bounds of approximation and the interval in which the corresponding inequality
holds true.
In this way, we get a class of inequalities as well as intervals in which these inequalities are true.
We believe that our method can  contribute to  broader application of the power series in the study of inequalities.

\bibstyle{plain}

\end{document}